\documentclass{icmart}

\contact[khonda@math.usc.edu]{University of Southern California, Los Angeles, CA 90089}

\numberwithin{equation}{section}

\newtheorem{theorem}{Theorem}[section]

\theoremstyle{definition}

\newtheorem{conj}[theorem]{Conjecture}
\newtheorem{condition}[theorem]{Condition}

\newcommand{\bdry}{\partial}

\newcommand{\Ker}{\operatorname{Ker}}

\title[The topology and geometry of contact structures in dimension three]{The topology and geometry of contact structures in dimension three}

\author[Ko Honda]{Ko Honda \thanks{Supported by an Alfred P.\ Sloan Fellowship and an NSF CAREER Award.}}

\begin{document}

\begin{abstract}
The goal of this article is to survey recent developments in the theory of contact structures in dimension three.
\end{abstract}



\maketitle

\vskip.4in

In this article we survey recent developments in three-dimensional contact geometry.  Three-dimensional contact geometry lies at the interface between 3- and 4-manifold geometries, and has been an essential part of the flurry in low-dimensional geometry and topology over the last 20 years.  In dimension 3, it relates to foliation theory and knot theory; in dimension 4, there are rich interactions with symplectic geometry.  In both dimensions, there are relations with gauge theories such as Seiberg-Witten theory and Heegaard Floer homology, as well as to dynamics.

\section{Tight vs. overtwisted}
A contact structure $\xi$ on a 3-manifold $M$ is a (maximally) nonintegrable 2-plane field distribution. In this paper we assume that $M$ is oriented and $\xi$ is the kernel of a global 1-form $\alpha$ which satisfies $\alpha\wedge d\alpha>0$.  (Such a contact structure is often called {\em coorientable}.)  Although $\xi$ is locally $\Ker(dz-ydx)$ by a classical theorem of Pfaff-Darboux and hence has {\em no local geometry}, the global study of contact structures is rather complicated, in a way that echoes the intricacies of symplectic geometry.

One of the fundamental questions is to determine $\pi_0$ of the space $Cont(M)$ of contact 2-plane fields on $M$ --- this is often called the ``classification'' of contact structures on $M$.  The work of Bennequin~\cite{Be} (later clarified by Eliashberg~\cite{El1}) indicated that contact structures, in dimension three, come in two flavors: {\em tight} and {\em overtwisted}.  We define an {\em overtwisted disk} to be an embedded disk $D\subset M$ such that $\xi_x=T_xD$ for all $x\in \bdry D$.  An overtwisted contact structure is one which admits an {\em overtwisted disk}, whereas a tight contact structure is one which does not.  What Bennequin showed is that the local model $(\mathbb{R}^3, \Ker(dz-ydx))$ is tight --- hence locally every contact structure is tight, although globally it may not be.   (Showing that a contact structure is tight is highly nontrivial, because one must show that {\em no} overtwisted disk exists, no matter how complicated the embedding!) Let $Cont^{OT}(M)$ be the space of overtwisted contact 2-plane fields on $M$ and $Cont^{Tight}(M)$ be the space of tight contact 2-plane fields on $M$.  Eliashberg showed in \cite{El2} that $\pi_0(Cont^{OT}(M))$ is the same as the homotopy classes of 2-plane fields on $M$.

The space of tight contact structures, on the other hand, is more intimately related to the topology of $M$.  Eliashberg \cite{El1} gave the first classification result for tight contact structures, namely that $Cont^{Tight}(S^3)$ is connected.  The analysis of tight contact structures on various 3-manifolds has become more manageable in recent years, with the introduction of {\em convex surfaces} by Giroux~\cite{Gi1} and {\em bypasses} by the author~\cite{H1}.  The world of tight contact structures, as we understand it now, is a veritable zoo!

Two important subcategories of tight contact structures are the {\em weakly symplectically fillable} ones and the {\em Stein fillable} ones.  In the former case, $(M,\xi)$ bounds a symplectic 4-manifold $(X,\omega)$ and $\omega|_{\xi}>0$.  In the latter, $(M,\xi)$ bounds a Stein domain $(X,\omega,J)$ and $\omega=d\alpha$ on $M$ for a contact 1-form $\alpha$ that defines $\xi$.  Fillable contact structures (of either type) are tight by a theorem of Gromov~\cite{Gr} and Eliashberg~\cite{El3}; this was proved using Gromov's theory of $J$-holomorphic curves.  Prototypical examples of weakly symplectically fillable contact structures are the perturbations of taut codimension 1 foliations, as explained in Eliashberg-Thurston~\cite{ET}.  Etnyre and author~\cite{EH2} showed that there exist tight contact structures which are not weakly symplectically fillable.  Other examples were later obtained by Lisca-Stipsicz~\cite{LS1,LS2}. Eliashberg~\cite{El4} showed that there are weakly symplectically fillable contact structures on the 3-torus $T^3$ which are not Stein fillable.  Further examples were given on torus bundles by Ding-Geiges~\cite{DG}.

$$
\begin{array}{|ccccc|}
\hline
& & & &\\
\mbox{ Tight} & \supsetneqq & \mbox{Weakly symplectically fillable} & \supsetneqq & \mbox{Stein fillable }\\
& & & &\\
\hline
\end{array}
$$

It is known that not every 3-manifold admits a tight contact structure. Etnyre and the author~\cite{EH1} showed that the Poincar\'e homology sphere with orientation opposite to the one induced on the link of an algebraic singularity has no tight contact structure.  Lisca and Stipsicz~\cite{LS3} have since shown that the Poincar\'e homology sphere can be incorporated into a larger class of small Seifert fibered spaces which do not admit tight contact structures.  Since all these examples of 3-manifolds without tight contact structures are Seifert fibered, it is natural to ask whether tight contact structures exist on all hyperbolic 3-manifolds.  It turns out that {\em universally tight} contact structures, i.e., contact structures $\xi$ on $M$ that pull back to tight contact structures on the universal cover of $M$, do not always exist on hyperbolic 3-manifolds~\cite{HKM3}.  Compare this to foliation theory where Roberts-Shareshian-Stein~\cite{RSS} have shown that there are infinitely many hyperbolic 3-manifolds which do not admit Reebless codimension 1 foliations.  There is related work of Calegari-Dunfield~\cite{CD} and Fenley~\cite{Fe}, as well as a different approach using Seiberg-Witten Floer homology, due to Kronheimer-Mrowka-Ozsv\'ath-Szab\'o~\cite{KMOS}.  However, it is still conceivable that every hyperbolic 3-manifold has a tight contact structure.  (The Weeks manifold --- the closed hyperbolic 3-manifold with the smallest known volume --- {\em does} have Stein fillable contact structures with either orientation.  This can be easily seen by appealing to surgery on the Borromean rings as in Gompf~\cite{Go}.)

Next we turn our attention to the question of classification.  For simplicity, assume that $M$ is irreducible.  Colin~\cite{Co} and Kazez, Mati\'c and the author~\cite{HKM0}, independently, have shown that $\pi_0(Cont^{Tight}(M))$ is infinite for a {\em toroidal} 3-manifold $M$, namely one which admits an embedded torus for which $\pi_1(T) \hookrightarrow \pi_1(M)$.  On the other hand, if $M$ is not toroidal, then $\pi_0(Cont^{Tight}(M))$ is finite by a theorem of Colin, Giroux and the author~\cite{CGH1,CGH2}.  The latter generalizes an earlier theorem, due to Kronheimer-Mrowka~\cite{KM0}, which states that there are finitely many homotopy classes of 2-plane fields which carry symplectically fillable contact structures.

\section{Open book decompositions}

A fundamental advance in contact geometry is the work of Giroux~\cite{Gi2} (building on earlier work of Thurston-Winkeln\-kemper \cite{TW}, Bennequin \cite{Be}, Eliashberg-Gromov \cite{EG}, and Torisu \cite{To}), which relates contact structures and open book decompositions.  We briefly summarize this work, and then describe the developments that have taken place since Giroux's ICM 2002 article~\cite{Gi2}.

Let $(S,h)$ be a pair consisting of a compact oriented surface $S$ with nonempty boundary and a diffeomorphism $h:S \rightarrow S$ which restricts to the identity on $\bdry S$, and let $K$ be a link in a closed oriented 3-manifold $M$. An {\em open book decomposition} for $M$ with {\em binding} $K$ is a homeomorphism between $((S \times[0,1]) /_{\sim _h}, (\bdry S \times[0,1]) /_{\sim_h})$ and $(M,K)$. The equivalence relation $\sim_h$ is generated by $(x,1) \sim_h (h(x),0)$ for $x\in S$ and $(y,t) \sim_h (y,t')$ for $y \in \bdry S$. We will often identify $M$ with $(S \times[0,1]) / _{\sim _h}$; with this identification $S_t= S \times \{t\}, t\in [0,1]$, is called a {\em page} of the open book decomposition and $h$ is called the {\em monodromy map}. Two open book decompositions are {\em equivalent} if there is an ambient isotopy taking binding to binding and pages to pages. We will denote an open book decomposition by $(S,h)$, although, strictly speaking, an open book decomposition is determined by the triple $(S,h,K)$. There is a slight difference --- if we do not specify $K\subset M$, we are referring to isomorphism classes of open books instead of isotopy classes.

Every closed 3-manifold has an open book decomposition, but it is not unique. One way of obtaining a different open book decomposition of the same manifold is to perform a positive stabilization. $(S',h')$ is a {\em positive stabilization} of $(S,h)$ if $S'$ is the union of the surface $S$ and a band $B$ attached along the boundary of $S$ (i.e., $S'$ is obtained from $S$ by attaching a 1-handle along $\bdry S$), and $h'$ is defined as follows. Let $\gamma$ be a simple closed curve in $S'$ ``dual'' to the cocore of $B$ (i.e., $\gamma$ intersects the cocore of $B$ at exactly one point) and let $id_B \cup h$ be the extension of $h$ by the identity map to $B \cup S$. Also let $R_\gamma$ be a {\em positive} or {\em right-handed} Dehn twist about $\gamma$.  Then for a {\em positive} stabilization $h'$ is given by $R_\gamma \circ (id_B \cup h)$. It is well-known that, if $(S',h')$ is a positive stabilization of an open book decomposition $(S,h)$ of $(M,K)$, then $(S',h')$ is an open book decomposition of $(M, K')$ where $K'$ is obtained by a Murasugi sum of $K$ (also called the {\em plumbing} of $K$) with a positive Hopf link.

A contact structure $\xi$ is said to be {\em supported} by the open book decomposition $(S,h,K)$ if there is a contact 1-form $\alpha$ satisfying the following:
\begin{enumerate}
\item $d\alpha$ restricts to a symplectic form on each fiber $S_t$;
\item $K$ is transverse to $\xi$, and the orientation on $K$ given by $\alpha$ is the same as the boundary orientation induced from $S$ coming from the symplectic structure.
\end{enumerate}
Thurston and Winkelnkemper \cite{TW} showed that any open book decomposition $(S,h,K)$ of $M$ supports a contact structure $\xi$. Moreover, the contact planes can be made arbitrarily close to the tangent planes of the pages (away from the binding).

The following result provides a converse (and more), due to Giroux~\cite{Gi2}.

\begin{theorem}[Giroux]
Every contact structure $(M,\xi)$ on a closed 3-manifold $M$ is supported by some open book decomposition $(S,h,K)$. Moreover, two open book decompositions $(S,h,K)$ and $(S',h',K')$ which support the same contact structure $(M,\xi)$ become equivalent after applying a sequence of positive stabilizations to each.
\end{theorem}

\subsection{Concave symplectic fillings}
Consider a closed 2-form $\omega_0$ on the contact 3-manifold $(M,\xi)$ for which $\omega_0|_\xi>0$. (Such a 2-form $\omega_0$ is often called a {\em dominating} 2-form for $\xi$.)  A {\em concave symplectic filling} for $(M,\xi,\omega_0)$ is a symplectic 4-manifold $(X,\omega)$ for which $\bdry X=-M$ and $i^*\omega=\omega_0$, where $i:M\rightarrow X$ is the inclusion.

The use of open book decompositions enabled Eliashberg~\cite{El5} and Etnyre~\cite{Et1} to construct concave symplectic fillings for any contact 3-manifold $(M,\xi)$ together with a dominating 2-form $\omega_0$. This turned out to be the only missing ingredient in Kronheimer-Mrowka's proof of Property P for knots~\cite{KM}.

\begin{theorem}[Kronheimer-Mrowka]
If $K\subset S^3$ is a nontrivial knot and $S^3_1(K)$ is the three-manifold obtained by $+1$-surgery along $K$, then $\pi_1(S^3_1(K))\not=0$.
\end{theorem}

In the 1980's Gabai~\cite{Ga} proved that $M=S^3_0(K)$ admits a taut foliation $\mathcal{F}$ ($\not=$ the foliation of $S^1\times S^2$ by $\{pt\}\times S^2$) if $K$ is not the unknot.  By Eliashberg-Thurston~\cite{ET}, $\mathcal{F}$ can be perturbed into a pair $\xi_+$, $\xi_-$ of positive and negative contact structures, and $X=M\times[0,1]$ admits a symplectic structure $\omega$ for which $\omega|_{\xi_+}>0$ at $M\times\{1\}$ and $\omega|_{\xi_-}<0$ at $M\times\{0\}$.  This used to be called a {\em symplectic semi-filling} of $(M,\xi_+)$ since $\bdry X$ had more than one component. (We no longer have the need to use the ``semi'' terminology, thanks to Eliashberg and Etnyre.)  The work of Eliashberg and Etnyre enabled one to fill both of the components of $\bdry X$ so that $(M,\xi_+)$ was now embedded in a closed symplectic manifold $X'$. Kronheimer and Mrowka were then able to appeal to: (i) the work by Taubes~\cite{Ta} on the nontriviality of Seiberg-Witten invariants of $X'$; (ii) the work by Feehan and Leness (see~\cite{FL}, for example) relating the Seiberg-Witten and Donaldson invariants; (iii) a stretching argument in instanton Floer homology; and (iv) Floer's exact triangle~\cite{Fl} for instanton Floer homology.

Another application of the existence of concave symplectic fillings is progress by Etnyre~\cite{Et2} on the following problem:  Given a contact manifold $(M,\xi)$, what is the minimum genus amongst all the pages of open books corresponding to $\xi$?  Etnyre has shown that many interesting classes of tight contact structures (among them perturbations of taut foliations) do not admit planar open book decompositions.

\subsection{Heegaard Floer homology}
Another important application of the open book framework is the definition of the {\em contact class} $c(\xi)$ in the Heegaard Floer homology of Ozsv\'ath-Szab\'o \cite{OSz1,OSz2}. Using open book decompositions, Ozsv\'ath and Szab\'o \cite{OSz3} defined an invariant of the contact structure $(M,\xi)$, which is an element $c(\xi)\in \widehat{HF}(M)$.  Among the many properties enjoyed by $c(\xi)$, we have the following:
\begin{enumerate}
\item If $\xi$ is overtwisted, then $c(\xi)=0$.
\item If $(M',\xi')$ is obtained from $(M,\xi)$ by Legendrian ($-1$) surgery, and $c(\xi)\not=0$, then $c(\xi')\not=0$.
\item If $\xi$ is weakly symplectically fillable, then $c(\xi)\not=0$ (provided ``twisted'' coefficients are used).
\end{enumerate}

Lisca and Stipsicz~\cite{LS4} showed that the contact class was surprisingly good at detecting tight contact structures --- Heegaard Floer homology could now be used to prove the tightness of many contact structures which were hitherto only conjectured to be tight.  This area is currently an active area of research, with contributions from Ghiggini~\cite{Gh1,Gh2}, Plamenevskaya~\cite{Pl1,Pl2}, etc.

\section{Right-veering}

We will now seek to explain the roles of tightness, weak symplectic fillability, and Stein fillability in the open book context.  Except for Theorem~\ref{LPG}, this is joint work with W.\ Kazez and G.\ Mati\'c and further details can be found in \cite{HKM1,HKM2}.

Let $S$ be a compact oriented surface with nonempty boundary. Denote by $Aut(S,\bdry S)$ the group of (isotopy classes of) diffeomorphisms of $S$ which restrict to the identity on $\bdry S$.  We have the monoid $Dehn^+(S,\bdry S)\subset Aut(S,\bdry S)$ of products of positive Dehn twists.  The following is due to Giroux~\cite{Gi2}, inspired by the work of Loi-Piergallini~\cite{LP}:

\begin{theorem}[Giroux]\label{LPG}
A contact structure $\xi$ on $M$ is Stein fillable if and only if $\xi$ is supported by some open book $(S,h,K)$ with $h\in Dehn^+(S,\bdry S)$.
\end{theorem}

We remark that the theorem does not say that {\em every} open book $(S,h)$ for $(M,\xi)$ Stein fillable satisfies $h\in Dehn^+(S,\bdry S)$.

There is another monoid, namely the monoid $Veer(S,\bdry S)$ of {\em right-veering} diffeomorphisms, which is intimately connected with the tight contact structures.  Given two properly embedded oriented arcs $\alpha$ and $\beta$ with the same initial point $x\in \bdry S$, we say $\alpha$ is {\em to the left} of $\beta$ if the following holds: Isotop $\alpha$ and $\beta$, while fixing their endpoints, so that they intersect transversely (this include the endpoints) and with the fewest possible number of intersections.  We then say $\alpha$ is to the left of $\beta$ if either $\alpha=\beta$ or the tangent vectors $(\dot \beta(0), \dot\alpha(0))$ define the orientation on $S$ at $x$.  Then $h$ is {\em right-veering} if for every choice of basepoint $x \in \bdry S$ and every choice of $\alpha$ based at $x$, $h(\alpha)$ is to the right of $\alpha$.  One easily sees that $Dehn^+(S,\bdry S)\subset Veer(S,\bdry S)$.

\begin{theorem}[Honda-Kazez-Mati\'c~\cite{HKM1}] \label{thm:veer}
A contact structure $(M,\xi)$ is tight if and only if all of its open book decompositions $(S,h)$ are such that $h\in Veer(S,\bdry S)$.
\end{theorem}

This theorem is an improvement over the ``sobering arc'' criterion for overtwistedness, given by Goodman~\cite{Goo}.

Now recall Thurston's classification of surface diffeomorphisms~\cite{Th}, which improved upon earlier work of Nielsen \cite{Ni1, Ni2, Ni3}. A diffeomorphism $h:S\rightarrow S$ satisfies one of the following:

\begin{enumerate}
\item $h$ is {\em reducible}, i.e., there exists an essential multicurve $\gamma$ such that $h(\gamma)$ is isotopic to $\gamma$.
\item $h$ is homotopic to a {\em periodic} homeomorphism $\psi$, i.e., there is an integer $n>0$ such that $\psi^n=id$.
\item $h$ is homotopic to a {\em pseudo-Anosov} homeomorphism $\psi$.
\end{enumerate}

We will now define the {\em fractional Dehn twist coefficients}, extensively studied by Gabai and Oertel (see for example \cite{GO}) in the context of essential laminations. Suppose for simplicity that $\bdry S$ is connected and $h\in Aut(S,\bdry S)$ is homotopic to a pseudo-Anosov representative $\psi$.  (The periodic case is analogous.)  Let $H: S\times[0,1]\rightarrow S$ be an isotopy from $h$ to $\psi$, i.e., $H(x,0)=h(x)$ and $H(x,1)=\psi(x)$. On the boundary of $S$, $\psi$ has $2n$ fixed points, $n$ attracting and $n$ repelling.  Let us label the attracting fixed points $x_1,\dots,x_n$ in order around $\bdry S$.  Now define $\beta: \bdry S \times [0,1] \to \bdry S\times [0,1]$ by sending $(x,t)\mapsto (H(x,t),t)$.  It follows that the arc $\beta(x_i \times [0,1])$ connects $(x_i, 0)$ and $(x_{i+k}, 1)$, for some $k$. We call $\beta$ a {\it fractional Dehn twist} by an amount $c\in \mathbb{Q}$, where $c\equiv k/n$ modulo $1$ is the number of times $\beta(x_i \times [0,1])$ circles around $\bdry S \times [0,1]$ (here circling in the direction of $\bdry C$ is considered positive). Form the union of $\bdry S \times [0,1]$ and $S$ by gluing $\bdry S \times \{1\}$ and $\bdry S$. By identifying this union with $S$, we construct the homeomorphism $\beta\cup\psi$ on $S$ which is isotopic to $h$, relative to $\bdry S$. (We will assume that $h=\beta\cup\psi$, although $\psi$ is usually just a homeomorphism, not a diffeomorphism.)

\begin{theorem}[Honda-Kazez-Mati\'c~\cite{HKM1}]\label{periodic-pA}
If $h$ is periodic, then $h$ is right-veering if and only if $c\geq 0$.  If $h$ is pseudo-Anosov, then $h$ is right-veering if and only if $c>0$.
\end{theorem}

Theorem~\ref{thm:veer} is not completely satisfactory --- ideally one should just need to look at one $(S,h)$ (instead of its equivalence class under stabilizations) to determine whether $(S,h)$ is tight, fillable, etc.  To this end, let us consider the case of the once-punctured torus $S$.  Suppose $h$ is pseudo-Anosov.  Then the following hold:

\begin{enumerate}
\item If $c\leq 0$, then $h$ is overtwisted.
\item If $c> 0$, then $h$ is tight.
\item If $c\geq 1$, then $h$ is weakly symplectically fillable and universally tight.
\item For any $c>0$ there exist $h\in Veer(S,\bdry S)-Dehn^+(S,\bdry S)$ whose fractional Dehn twist coefficient is equal to $c$.
\end{enumerate}

For the once-punctured torus, $c>0$ is equivalent to $c\geq {1\over 2}$, since the pseudo-Anosov representative $\psi$ will have $n=2$ attracting fixed points.  (2) is proved using Heegaard-Floer homology.  (3) is proved by showing that the taut foliations constructed by Roberts in \cite{Ro1,Ro2} can be perturbed to the contact structure adapted to the open book. (4) is proved by looking at a function on the Farey tessellation called the {\em Rademacher function} (see for example~\cite{GG}).  It is very plausible that many, if not all, of the $h\in Veer(S,\bdry S)-Dehn^+(S,\bdry S)$ never become products of positive Dehn twists after (repeated) stabilization, and that such $(S,h)$ are indeed not Stein fillable.

Once we restrict our attention to right-veering $(S,h)$, calculations in {\em contact homology} and {\em Heegaard Floer homology} both become more manageable.  In the rest of the paper, we focus on contact homology, leaving the Heegaard Floer aspects for another occasion.

\section{Contact homology}

Given a contact form $\alpha$ for $(M,\xi)$, there is a corresponding {\em Reeb vector field} $R$ defined as follows: $i_{R} d\alpha=0$, $i_{R}\alpha=1$.  One of the motivating questions in the study of the dynamics of Reeb vector fields is the following Weinstein conjecture (in dimension three):

\begin{conj}[Weinstein conjecture]
Let $(M,\xi)$ be a contact 3-manifold.  Then for any contact form $\alpha$ with $\Ker(\alpha)=\xi$, the corresponding Reeb vector field $R$ admits a closed periodic orbit.
\end{conj}

The fundamental step was taken when Hofer~\cite{Ho1} studied $J$-holomorphic disks in the symplectization $(\mathbb{R}\times M, d(e^t\alpha))$, and proved that there is always bubbling (and hence a closed periodic orbit) when $(M,\xi)$ is overtwisted or $\pi_2(M)\not=0$.  This showed that the Weinstein conjecture holds for overtwisted contact structures and 3-manifolds with $\pi_2(M)\not=0$.  (Hofer also showed that the Weinstein conjecture holds for $S^3$.)  Hofer's work has subsequently bubbled off a large industry in contact dynamics, and we mention only a few highlights.  The properties of holomorphic curves in symplectizations were analyzed by Hofer-Wysocki-Zehnder~\cite{HWZ1,HWZ2,HWZ3}.  Also Eliashberg, Givental and Hofer~\cite{EGH} have suggested a {\em Symplectic Field Theory}, a Floer-type theory involving closed orbits of Reeb vector fields and holomorphic curves ``bounding'' these closed orbits.  The technical details of this theory have finally started to appear --- see~\cite{HWZ4}.

The Weinstein conjecture in dimension three has been verified for contact structures which admit planar open book decompositions~\cite{ACH} (also see related work of Etnyre~\cite{Et2}), for certain Stein fillable contact structures~\cite{Ch,Ze}, and for certain universally tight contact structures on toroidal manifolds~\cite{BC}.  We also refer the reader to the survey article by Hofer~\cite{Ho2}.

In \cite{CH1}, Colin and the author prove the following:

\begin{theorem}[Colin-Honda~\cite{CH1}]
The Weinstein conjecture holds for contact structures $(M,\xi)$ which have open books $(S,h)$ with periodic monodromy. (Here $S$ may have many boundary components.)
\end{theorem}

Suppose $\bdry S$ is connected.  If the fractional Dehn twist coefficient $c<0$, then $(M,\xi)$ is overtwisted by Theorem~\ref{periodic-pA} above.  If $c=0$, then $M$ is a connected sum, and if $0<c< {1\over 2g(S)-1}$, then the universal cover of $M$ is $S^3$.  Here $g(S)$ is the genus of the closed surface obtained by capping off $S$ with a disk.  In all the above cases, the Weinstein conjecture has been settled by Hofer~\cite{Ho1}.

It remains to examine the case where $c\geq {1\over 2g(S)-1}$.  The following is the main result in \cite{CH1}:

\begin{theorem}[Colin-Honda~\cite{CH1}]\label{nontrivial}
The cylindrical contact homology of $(M,\xi)$, as defined below, exists and is nonzero if $c\geq {1\over 2g(S)-1}$.
\end{theorem}

The nontriviality of the cylindrical contact homology for $(M,\xi)$ implies the Weinstein conjecture for $(M,\xi)$.

{\em Contact homology} is the simplest version of Symplectic Field Theory which takes place in the symplectization $(\mathbb{R}\times M,d(e^t\alpha))$ and counts punctured $J$-holo\-morphic spheres $\tilde u:\Sigma\rightarrow \mathbb{R}\times M$ with one positive end and many negative ends. (Here $t$ is the coordinate for $\mathbb{R}$.)  {\em Cylindrical} contact homology, if it exists, is a version which only counts $J$-holomorphic cylinders, i.e., spheres with two punctures. The cylindrical theory exists, if the following condition holds:

\begin{condition}\label{cylindrical}
There exists a nondegenerate Reeb vector field $R$ for which no contractible periodic orbit $\gamma$ with Conley-Zehnder index $\mu(\gamma)=0$, $1$ or $2$ bounds a finite energy plane (at the positive end) in the symplectization $\mathbb{R}\times M$.
\end{condition}

Assuming Condition~\ref{cylindrical}, we now define the {\em cylindrical contact homology}.

Let $\alpha$ be a contact 1-form for which $R$ is a nondegenerate Reeb vector field, and let $J$ be an almost complex structure on $\mathbb{R}\times M$ which is {\em adapted} to the symplectization: If we write $T_{(t,x)} (\mathbb{R}\times M)=\mathbb{R}{\bdry\over \bdry t} \oplus \mathbb{R} R\oplus \xi$, then $J$ maps $\xi$ to itself and sends ${\bdry\over \bdry t}\mapsto R$, $R\mapsto -{\bdry\over \bdry t}$.  Let $\mathcal{P}$ be the collection of closed orbits of $R$. (We may need to omit certain closed orbits, but we will not worry about these technicalities here.)

If $\gamma$ is a contractible periodic orbit which bounds a disk $D$, then we trivialize $\xi|_D$ and define the {\em Conley-Zehnder index} $\mu(\gamma,D)$ to be the Conley-Zehnder index of the path of symplectic maps $\{d\phi_t: \xi_{\gamma(0)}\rightarrow \xi_{\gamma(t)}, t\in[0,T]\}$ with respect to this trivialization, where $\phi_t$ is the time $t$ flow of the Reeb vector field $R$ and $T$ is the period of $\gamma$. In our case, $M$ is Seifert fibered and $\pi_2(M)=0$, so $\mu(\gamma)$ is independent of the choice of $D$.  If $\gamma$, $\gamma'$ are not contractible, but belong to the same {\em free homotopy class} $[\gamma]=[\gamma']$, then let $Z$ be the cylinder between $\gamma$ and $\gamma'$.  Trivialize $\xi|_Z$ and define the {\em relative} Conley-Zehnder index $\mu(\gamma,\gamma')$ to be the Conley-Zehnder index of $\gamma$ minus the Conley-Zehnder index of $\gamma'$, both calculated with respect to this trivialization.  Again, in our case, $\mu(\gamma,\gamma')$ does not depend on our choice of $Z$.

Define the moduli space
$$\mathcal{M}(J,\gamma_+,\gamma_-)=\left\{
\begin{array}{c}
\mbox{$J$-holomorphic cylinders $\tilde u=(a,u):\mathbb{R}\times S^1\rightarrow \mathbb{R}\times M$}\\
\lim_{s\rightarrow \pm \infty} u(s,t) = \gamma_\pm(t)\\
\lim_{s\rightarrow \pm \infty} a(s,t) = \pm \infty
\end{array}
\right\}.$$
Here, $\gamma_\pm(t)$ refers to some parametrization of the trajectory $\gamma_\pm$.  The convergence for $u(s,t)$ and $a(s,t)$ is in the $C^0$-topology.  The complex structure $j$ on $\mathbb{R}\times S^1$ is the usual one: if $(s,t)$ are coordinates on $\mathbb{R}\times \mathbb{R}/\mathbb{Z}$, then $j:{\bdry \over \bdry s}\mapsto {\bdry\over \bdry t}$.  We choose a {\em regular} $J$ (still adapted to the symplectization) for which $\mathcal{M}(J,\gamma_+,\gamma_-)$ is a transverse zero set of the $\overline\bdry$-operator and has the expected dimension $\mu(\gamma_+,\gamma_-)$.

The chain group is the $\mathbb{Q}$-vector space $C=\mathbb{Q}\langle \mathcal{P} \rangle$ generated by $\mathcal{P}$.  Now the boundary map $\bdry: C\rightarrow C$ is given on elements $\gamma\in\mathcal{P}$ by:
$$\bdry\gamma= \sum_{\gamma'\in\mathcal{P}, ~\mu(\gamma,\gamma')=1} {n_{\gamma,\gamma'}\over \kappa(\gamma')} ~~\gamma',
$$
where $\kappa(\gamma)$ is the multiplicity of $\gamma$. If $\mu(\gamma,\gamma')=1$, then $\mathcal{M}(\gamma,\gamma')$ is a 1-dimensional moduli space and we quotient out by translations in the $\mathbb{R}$-direction.  Then $n_{\gamma,\gamma'}$ is a signed count of points in $\mathcal{M}(\gamma,\gamma')/\mathbb{R}$, following a coherent orientation scheme given in \cite{EGH}.  If $\gamma$, $\gamma'$ are multiply covered, then each non-multiply-covered holomorphic curve $\tilde u\in \mathcal{M}(\gamma,\gamma')/\mathbb{R}$ contributes $\pm \kappa(\gamma) \kappa(\gamma')$ to $n_{\gamma,\gamma'}$. If $\tilde u$ is a $k$-fold cover of a somewhere injective holomorphic curve, then it is counted as $\pm k(\kappa(\gamma) \kappa(\gamma'))$.  The definition of $\bdry$ is extended linearly to all of $C$.  (For the purposes of Theorem~\ref{nontrivial}, we can restrict attention to the portion of $\mathcal{P}$ consisting of non-multiply-covered orbits, so we may work with $\mathbb{Z}/2\mathbb{Z}$-coefficients.)

With the above restriction (Condition~\ref{cylindrical}) on the Conley-Zehnder indices of contractible periodic orbits, it can be shown that $\bdry\circ\bdry=0$, hence the cylindrical contact homology is well-defined.  Moreover, it does not depend on the choice of generic $J$ or on the choice of nondegenerate $R$, provided Condition~\ref{cylindrical} is satisfied.

\vskip.2in
We now indicate some elements of the proof of Theorem~\ref{nontrivial}.  The well-definition of cylindrical contact homology is proved by projecting a finite energy plane $\tilde u =(a,u):\mathbb{R}^2 \rightarrow \mathbb{R}\times M$ to $M$, observing that $u:\mathbb{R}^2\rightarrow M$ is positively transverse to the Reeb vector field $R$ except at complex branch points, and using (a generalization of) the Rademacher function.  A similar technique gives restrictions on pairs $\gamma$, $\gamma'$ which admit holomorphic cylinders between them.  Since we have enough restrictions on the boundary maps, an Euler characteristic argument gives the result.

In the pseudo-Anosov case we expect the analog of Theorem~\ref{nontrivial} to hold when $c>{1\over n}$, where $n$ is the number of attracting (= number of repelling) periodic points. (These are currently being worked out in \cite{CH2}.)  It still remains to consider $c={1\over n}$ in the pseudo-Anosov case....

\vskip.2in
\noindent
{\em Acknowledgements:}  The author wholeheartedly thanks Francis Bonahon for many suggestions on improving the exposition.


\frenchspacing
\begin{thebibliography}{7}

\bibitem{ACH}
C.\ Abbas, K.\ Cieliebak and H.\ Hofer, \textit{The Weinstein conjecture for planar contact structures in dimension three},  to appear in Comment.\ Math.\ Helv.

\bibitem{Be}
D.\ Bennequin, \textit{Entrelacements et \'equations de Pfaff}, Third Schnepfenried geometry conference, Vol.\ 1 (Schnepfenried, 1982), 87--161, Ast\'erisque \textbf{107--108}, Soc.\ Math.\ France, Paris,
1983.

\bibitem{BC}
F.\ Bourgeois and V.\ Colin, \textit{Homologie de contact des vari\'et\'es tor\"oidales}, Geom.\ Topol.\ {\bf 9} (2005) 299--313 (electronic).

\bibitem{CD}
D.\ Calegari and N.\ Dunfield, \textit{Laminations and groups of homeomorphisms of the circle}, Invent.\ Math.\ {\bf 152}  (2003), 149--204.

\bibitem{Ch}
W.\ Chen, \textit{Pseudo-holomorphic curves and the Weinstein conjecture},  Comm.\ Anal.\ Geom.\  {\bf 8}  (2000), 115--131.

\bibitem{Co}
V.\ Colin, \textit{Une infinit\'e de structures de contact tendues sur les vari\'et\'es toro\"\i dales}, Comment.\ Math.\ Helv.\ {\bf 76} (2001), 353--372.

\bibitem{CGH1}
V.\ Colin, E.\ Giroux and K.\ Honda, \textit{On the coarse classification of tight contact structures}, Topology and geometry of manifolds (Athens, GA, 2001),  109--120, Proc.\ Sympos.\ Pure Math., 71, Amer.\ Math.\ Soc., Providence, RI, 2003.

\bibitem{CGH2}
V.\ Colin, E.\ Giroux and K.\ Honda, \textit{Finitude homotopique et isotopique des structures de contact tendues}, in preparation.

\bibitem{CH1}
V.\ Colin and K.\ Honda, \textit{Reeb vector fields and open book decompositions I: the periodic case}, preprint 2005.

\bibitem{CH2}
V.\ Colin and K.\ Honda, \textit{Reeb vector fields and open book decompositions II: the pseudo-Anosov case}, in preparation.

\bibitem{DG}
F.\ Ding and H.\ Geiges, \textit{Symplectic fillability of tight contact structures on torus bundles},  Algebr.\ Geom.\ Topol.\  {\bf 1}  (2001), 153--172 (electronic).

\bibitem{El1}
Y.\ Eliashberg, \textit{Contact 3-manifolds twenty years since J.\ Martinet's work}, Ann.\ Inst.\ Fourier (Grenoble) \textbf{42} (1992), 165--192.

\bibitem{El2}
Y.\ Eliashberg, \textit{Classification of overtwisted contact structures on 3-manifolds}, Invent.\ Math.\ \textbf{98} (1989), 623--637.

\bibitem{El3}
Y.\ Eliashberg, \textit{Filling by holomorphic discs and its applications}, Geometry of low-dimensional
manifolds, 2 (Durham, 1989), 45--67, London Math.\ Soc.\ Lecture Note Ser.\
{\bf 151}, Cambridge Univ.\ Press, Cambridge, 1990.

\bibitem{El4}
Y.\ Eliashberg, \textit{Unique holomorphically fillable contact structure on the $3$-torus}, Internat.\ Math.\ Res.\ Notices {\bf 1996},  77--82.

\bibitem{El5}
Y.\ Eliashberg, \textit{A few remarks about symplectic filling},  Geom.\ Topol.\  {\bf 8}  (2004), 277--293 (electronic).

\bibitem{EGH}
Y.\ Eliashberg, A.\ Givental and H.\ Hofer,  {\it Introduction to symplectic field theory}, GAFA 2000 (Tel Aviv, 1999), Geom.\ Funct.\ Anal.\ 2000, Special Volume, Part II, 560--673.

\bibitem{EG}
Y.\ Eliashberg and M.\ Gromov, \textit{Convex symplectic manifolds}, Several complex variables and complex geometry, Part 2 (Santa Cruz, CA, 1989),  135--162, Proc.\ Sympos.\ Pure Math., 52, Part 2, Amer.\ Math.\ Soc., Providence, RI, 1991.

\bibitem{ET}
Y.\ Eliashberg and W.\ Thurston, \textit{Confoliations}, University Lecture Series \textbf{13}, Amer.\ Math.\ Soc., Providence, 1998.

\bibitem{Et1}
J.\ Etnyre, \textit{On symplectic fillings}, Algebr.\ Geom.\ Topol.\ {\bf 4}  (2004), 73--80 (electronic).

\bibitem{Et2}
J.\ Etnyre, \textit{Planar open book decompositions and contact structures},  Int.\ Math.\ Res.\ Not.\  {\bf 2004}, 4255--4267.

\bibitem{EH1}
J.\ Etnyre and K.\ Honda, \textit{On the nonexistence of tight contact structures}, Ann.\ of Math.\ (2) {\bf 153} (2001), 749--766.

\bibitem{EH2}
J.\ Etnyre and K.\ Honda, \textit{Tight contact structures with no symplectic fillings}, Invent.\ Math.\ {\bf 148} (2002), 609--626.

\bibitem{FL}
P.\ Feehan and T.\ Leness, \textit{On Donaldson and Seiberg-Witten invariants}, Topology and geometry of manifolds (Athens, GA, 2001),  237--248, Proc.\ Sympos.\ Pure Math., 71, Amer.\ Math.\ Soc., Providence, RI, 2003.

\bibitem{Fe}
S.\ Fenley, \textit{Laminar free hyperbolic 3-manifolds}, preprint 2002. \texttt{ArXiv:math.GT/0210482.}

\bibitem{Fl}
A.\ Floer, \textit{Instanton homology and Dehn surgery}, The Floer memorial volume, 77--97, Progr.\ Math., 133, Birkh\"auser, Basel, 1995.

\bibitem{Ga}
D.\ Gabai, \textit{Foliations and the topology of 3-manifolds. II}, J.\ Differential Geom.\ {\bf 26} (1987), 461--478.

\bibitem{GO}
D.\ Gabai and U.\ Oertel, \textit{Essential laminations in $3$-manifolds},  Ann.\ of Math.\ (2) {\bf 130}  (1989),  41--73.

\bibitem{GG}
J.-M.\ Gambaudo and E.\ Ghys, \textit{Braids and signatures}, Bull.\ Soc.\ Math.\ France, to appear.

\bibitem{Gh1}
P.\ Ghiggini, \textit{Strongly fillable contact 3-manifolds without Stein fillings}, Geom.\ Topol.\ {\bf 9}  (2005), 1677--1687 (electronic).

\bibitem{Gh2}
P.\ Ghiggini, \textit{Infinitely many universally tight contact manifolds with trivial Ozsv\'ath-Sz\'ab\'o contact invariants}, preprint 2005. \texttt{ArXiv:math.GT/0510574.}

\bibitem{Gi1}
E.\ Giroux, \textit{Convexit\'e en topologie de contact}, Comment.\ Math.\ Helv.\  {\bf 66} (1991), 637--677.

\bibitem{Gi2}
E.\ Giroux, \textit{G\'eom\'etrie de contact:\ de la dimension trois vers les dimensions sup\'erieures}, Proceedings of the International Congress of Mathematicians, Vol.\ II (Beijing, 2002), 405--414, Higher Ed.\ Press, Beijing, 2002.

\bibitem{Go}
R.\ Gompf, \textit{Handlebody construction of Stein surfaces}, Ann.\ of Math.\ (2) \textbf{148} (1998), 619--693.

\bibitem{Goo}
N.\ Goodman, \textit{Overtwisted open books from sobering arcs}, Algebr.\ Geom.\ Topol.\ {\bf 5} (2005), 1173--1195 (electronic).

\bibitem{Gr}
M.\ Gromov, \textit{Pseudo-holomorphic curves in symplectic manifolds}, Invent.\ Math.\ {\bf 82} (1985), 307--347.

\bibitem{Ho1}
H.\ Hofer, \textit{Pseudoholomorphic curves in symplectizations with applications to the Weinstein conjecture in dimension three}, Invent.\ Math.\ {\bf 114} (1993), 515--563.

\bibitem{Ho2}
H.\ Hofer, \textit{Dynamics, topology, and holomorphic curves}, Proceedings of the International Congress of Mathematicians, Vol.\ I (Berlin, 1998).  Doc.\ Math.\  1998,  Extra Vol. I, 255--280 (electronic).

\bibitem{HWZ1}
H.\ Hofer, K.\ Wysocki and E.\ Zehnder, {\it Properties of pseudo-holomorphic curves in symplectizations I: asymptotics},  Ann.\ Inst.\ H.\ Poincar\'e Anal.\ Non Lin\'eaire  {\bf 13}  (1996),  337--379.

\bibitem{HWZ2}
H.\ Hofer, K.\ Wysocki and E.\ Zehnder, {\it Properties of pseudo-holomorphic curves in symplectisations II: embedding controls and algebraic invariants}, Geom.\ Funct.\ Anal.\  {\bf 5}  (1995), 270--328.

\bibitem{HWZ3}
H.\ Hofer, K.\ Wysocki and E.\ Zehnder, {\it Properties of pseudoholomorphic curves in symplectizations. III. Fredholm theory}, Topics in nonlinear analysis, 381--475, Progr.\ Nonlinear Differential Equations Appl.\ {\bf 35}, Birkh\"auser, Basel, 1999.

\bibitem{HWZ4}
H.\ Hofer, K.\ Wysocki and E.\ Zehnder, \textit{Polyfolds and Fredholm theory, Part I}, preprint 2005.

\bibitem{H1}
K.\ Honda, \textit{On the classification of tight contact structures I}, Geom.\ Topol.\ {\bf 4} (2000), 309--368 (electronic).

\bibitem{HKM0}
K.\ Honda, W.\ Kazez and G.\ Mati\'c, \textit{Convex decomposition theory}, Internat.\ Math.\ Res.\ Notices {\bf 2002}, 55--88.

\bibitem{HKM1}
K.\ Honda, W.\ Kazez and G.\ Mati\'c, \textit{Right-veering diffeomorphisms of compact surfaces with boundary I}, preprint 2005.

\bibitem{HKM2}
K.\ Honda, W.\ Kazez and G.\ Mati\'c, \textit{Right-veering diffeomorphisms of compact surfaces with boundary II}, preprint 2005.

\bibitem{HKM3}
K.\ Honda, W.\ Kazez and G.\ Mati\'c, in preparation.

\bibitem{KM0}
P.\ Kronheimer and T.\ Mrowka,  \textit{Monopoles and contact structures}, Invent.\ Math.\ \textbf{130} (1997), 209--255.

\bibitem{KM}
P.\ Kronheimer and T.\ Mrowka, \textit{Witten's conjecture and property P},  Geom.\ Topol.\ {\bf 8}  (2004), 295--310 (electronic).

\bibitem{KMOS}
P.\ Kronheimer, T.\ Mrowka, P.\ Ozsv\'ath and Z.\ Szab\'o, \textit{Monopoles and lens space surgeries}, to appear in Ann.\ of Math.

\bibitem{LS1}
P.\ Lisca and A.\ Stipsicz, \textit{An infinite family of tight, not semi-fillable contact three-manifolds}, Geom.\ Topol.\ {\bf 7}  (2003), 1055--1073 (electronic).

\bibitem{LS2}
P.\ Lisca and A.\ Stipsicz, \textit{Tight, not semi-fillable contact circle bundles}, Math.\ Ann.\ {\bf  328}  (2004), 285--298.

\bibitem{LS4}
P.\ Lisca and A.\ Stipsicz, \textit{Ozsv\'ath-Szab\'o invariants and tight contact three-manifolds. I}, Geom.\ Topol.\ {\bf 8}  (2004), 925--945 (electronic).

\bibitem{LS3}
P.\ Lisca and A.\ Stipsicz, \textit{Ozsv\'ath-Szab\'o invariants and tight contact three-manifolds, II}, preprint 2004. \texttt{ArXiv:math.SG/0404136.}

\bibitem{LP}
A.\ Loi and R.\ Piergallini, \textit{Compact Stein surfaces with boundary as branched covers of $B\sp 4$}, Invent.\ Math.\ {\bf 143} (2001), 325--348.

\bibitem{Ni1}
J.\ Nielsen, \textit{Untersuchungen zur Topologie der geschlossenen zweiseitigen Fl\"achen I}, Acta Math.\ {\bf 50} (1927), 189--358.

\bibitem{Ni2}
J.\ Nielsen, \textit{Untersuchungen zur Topologie der geschlossenen zweiseitigen Fl\"achen II}, Acta Math.\ {\bf 53} (1929), 1--76.

\bibitem{Ni3}
J.\ Nielsen, \textit{Untersuchungen zur Topologie der geschlossenen zweiseitigen Fl\"achen III}, Acta Math.\ {\bf 58} (1931), 87--167.

\bibitem{OSz1}
P.\ Ozsv\'ath and Z.\ Szab\'o, \textit{Holomorphic disks and topological invariants for closed three-manifolds},  Ann.\ of Math.\ (2) {\bf 159}  (2004),  1027--1158.

\bibitem{OSz2}
P.\ Ozsv\'ath and Z.\ Szab\'o, \textit{Holomorphic disks and three-manifold invariants: properties and applications}, Ann.\ of Math.\ (2) {\bf 159}  (2004),  1159--1245.

\bibitem{OSz3}
P.\ Ozsv\'ath and Z.\ Szab\'o, \textit{Heegaard Floer homology and contact structures},  Duke Math.\ J.\  129  (2005),  39--61.

\bibitem{Pl1}
O.\ Plamenevskaya, \textit{Contact structures with distinct Heegaard Floer invariants}, Math.\ Res.\ Lett.\ {\bf 11}  (2004), 547--561.

\bibitem{Pl2}
O.\ Plamenevskaya, \textit{Transverse knots, branched double covers and Heegaard Floer contact invariants}, preprint 2004. \texttt{ArXiv:math.GT/0412183}.

\bibitem{Ro1}
R.\ Roberts, \textit{Taut foliations in punctured surface bundles, I}, Proc.\ London Math.\ Soc.\ (3)  {\bf 82}  (2001), 747--768.

\bibitem{Ro2}
R.\ Roberts, \textit{Taut foliations in punctured surface bundles, II}, Proc.\ London Math.\ Soc.\ (3) {\bf 83} (2001), 443--471.

\bibitem{RSS}
R.\ Roberts, J.\ Shareshian and M.\ Stein, \textit{Infinitely many hyperbolic 3-manifolds which contain no Reebless foliation}, J.\ Amer.\ Math.\ Soc.\ {\bf 16}  (2003),  639--679.

\bibitem{Ta}
C.\ Taubes, \textit{The Seiberg-Witten invariants and symplectic forms}, Math.\ Res.\ Lett.\ {\bf 1}  (1994), 809--822.

\bibitem{Th}
W.\ Thurston, \textit{On the geometry and dynamics of diffeomorphisms of surfaces}, Bull.\ Amer.\ Math.\ Soc.\ {\bf 19} (1988), 417--431.

\bibitem{To}
I.\ Torisu, \textit{Convex contact structures and fibered links in 3-manifolds}, Int.\ Math.\ Res.\ Notices {\bf 2000}, 441--454.

\bibitem{TW}
W.\ Thurston and H.\ Winkelnkemper, \textit{On the existence of contact forms}, Proc.\ Amer.\ Math.\ Soc.\ {\bf 52} (1975), 345-- 347.

\bibitem{Ze}
K.\ Zehmisch, \textit{Strong fillability and the Weinstein conjecture}, ArXiv:math.SG/0405203.


\end{thebibliography}
\end{document}